\documentclass[12pt]{amsart}
\usepackage{amscd,amssymb}

\newtheorem{thm}{Theorem}[section]

\newtheorem{lemma}[thm]{Lemma}
\newtheorem{prop}[thm]{Proposition}
\theoremstyle{definition}

\theoremstyle{remark}

\numberwithin{equation}{section}
\theoremstyle{question}

\def\norm#1{\left\Vert#1\right\Vert}

\def\C {{\mathbf C}}

\def\Is{{\mathrm{Is}}\,}

\def\d{{\mathrm{dist}}}

\def\QED{\nobreak\quad\ifmmode\roman{Q.E.D.}\else{\rm Q.E.D.}\fi}

\def\d{\delta}

\def\R{{\mathbf R}}

\def\sbs{\subset}

\def\e{\epsilon}

\def\ti{\times}
\def\obr{^{-1}}
\def\stm{\setminus}

\begin{document}

\title
{On unitary representations of groups of isometries}

\author{V. V. Uspenskij}

\address{Department of Mathematics, 321 Morton Hall, Ohio
University, Athens, Ohio 45701, USA}

\email{uspensk@math.ohiou.edu}


\date{10 June 2004}


\begin{abstract}
Is the topological group of all motions (including translations) 
of an infinite-dimensional
Hilbert space $H$ isomorphic to a subgroup of the unitary group
$U(H)$? This question was asked by Su Gao. We answer the question
in the affirmative. 
\end{abstract}

\maketitle

\setcounter{tocdepth}{1}

\section{Introduction}

For a Hilbert space $H$ we denote by $U_s(H)$ the group of all unitary
operators on $H$, equipped with the strong operator topology. This topology
is the same as the topology of pointwise convergence, it is
induced on $U(H)$ by the product topology on $H^H$. The group $U_s(H)$
is a topological group, in other words, its topology is compatible with 
its group structure.

Given a topological group $G$, a {\em unitary representation} of $G$
is a continuous homomorphism $f:G\to U_s(H)$. The representation $f$ 
is {\em faithful} if it is injective, and {\em topologically faithful}
if $f$ is a homeomorphic embedding of $G$ into $U_s(H)$. Thus 
a group admits a topologically faithful unitary representation 
if and only if it is isomorphic to a topological subgroup of a group
of the form $U_s(H)$.

Every locally compact group $G$ (in particular, every Lie group)
admits a topologically faithful unitary representation. For
example, the regular representation on $L^2(G)$ is such. Beyond the class
of locally compact groups the situation is less clear. There
are topological groups with a countable base which do not admit 
any faithful unitary representation or even any faithful representation
by isometries on a reflexive Banach space (Megrelishvili). 
For example, the group
$H(I)$ of all self-homeomorphisms of the interval $I=[0,1]$ has this
property \cite{Megr}. It follows that $\Is(U)$, the group of isometries 
of the Urysohn universal metric space $U$, also has this property, 
since it contains a topological copy of $H(I)$. See \cite{UspCMUC}, 
\cite{UspArchive}, \cite{UspPr} for the proof and a discussion of 
the space $U$ and its group of isometries. 

The question arises: what are metric spaces $M$ for which the group
$\Is(M)$ of all isometries of $M$ admits a (topologically) faithful
unitary representation? We consider the topology of pointwise convergence
on $\Is(M)$ (or the compact-open topology, which is the same on $\Is(M)$). 
If $M$ is locally compact and connected, the group $\Is(M)$ is locally compact
\cite[Ch.1, Thm.4.7]{KN} and hence embeds in a unitary group. On the other
hand, if $M$ is the Urysohn space $U$ (which topologically is the same
as a Hilbert space \cite{UspUr}), then $\Is(M)$ has no faithful unitary
representations, as we noted above.

Let $H$ be a (real or complex) Hilbert space. Consider the topological
group $\Is(H)$ of all (not necessarily linear) isometries of $H$.
This group is a topological semidirect product of the subgroup of translations
(which is isomorphic to the additive group of $H$) and the unitary group
$U(H_\R)$ of the real Hilbert space $H_\R$ underlying $H$. 
Su Gao asked the author whether the topological group $\Is(H)$ is isomorphic
to a subgroup of a unitary group. The aim of the present 
note is to answer this question in the affirmative. 
It follows that the additive group of $H$ admits a topologically 
faithful unitary representation. Even this corollary does not look self-evident.

\begin{thm}
\label{main}
Let $H$ be a Hilbert space, and let $G=\Is(H)$ be the topological
group of all (not necessarily linear) isometries of $H$. Then $G$
is isomorphic to a subgroup of a unitary group.
\end{thm}

In Section~\ref{2} we remind some basic facts concerning positive-definite
functions on groups, and prove (a generalization of) Theorem~\ref{main}
in Section~\ref{3}.

\section{Positive-definite functions}
\label{2}
For a complex matrix $A=(a_{ij})$ we denote by $A^*$ the matrix $(\bar a_{ji})$.
$A$ is {\em Hermitian} if $A=A^*$. A Hermitian matrix $A$ is {\em positive}
%(notation: $A\ge 0$) 
if all the eigenvalues of $A$ are $\ge 0$ or, equivalently,
if $A=B^2$ for some Hermitian $B$. A complex function $p$ on a group $G$ is
{\em positive-definite} if for every $g_1,\dots, g_n\in G$ the $n\ti n$-matrix
$p(g_i\obr g_j)$ is Hermitian and positive. If $f:G\to U(H)$ is a unitary
representation of $G$, then for every vector $v\in H$ the function $p=p_v$
on $G$ defined by $p(g)=(gv,v)$ is positive-definite. 
(We denote by $(x,y)$ the scalar product of $x,y\in H$.) 
Conversely, let $p$ be a positive-definite function on $G$.
Then $p=p_v$ for some unitary representation $f:G\to U(H)$ and some $v\in H$.
Indeed, consider the group algebra $\C[G]$, equip it with the scalar
product defined by $(g,h)=p(h\obr g)$, quotient out the kernel, and take
the completion to get $H$. The regular representation of $G$ on $\C[G]$
gives rise to a unitary representation on $H$ with the required property.
If $G$ is a topological group and the function $p$ is continuous, the resulting
unitary representation is continuous as well. These arguments yield the
following criterion:

\begin{prop}
\label{crit}
A topological group $G$ admits a topologically faithful unitary representation
(in other words, is isomorphic to a subgroup of $U_s(H)$ for some Hilbert space
$H$) if and only if for every neighborhood $U$ of the neutral element $e$ of $G$
there exist a continuous positive-definite function $p:G\to \C$ and $a>0$
such that $p(e)=1$ and $|1-p(g)|>a$ for every $g\in G\stm U$.
\end{prop}

\begin{proof} The necessity is easy: consider convex linear combinations
of the functions $p_v$ defined above, where $v\in H$ is a unit vector. We prove
that the condition is sufficient. According to the remarks preceding the
proposition, for every neighborhood $U$ of $e$ there exists 
a unitary representation $f_U: G\to U_s(H)$ such that $f_U\obr(V)\sbs U$
for some neighborhood $V$ of $1$ in $U_s(H)$. Indeed, pick a positive-definite
function $p$ on $G$ such that $p(e)=1$ and 
$$
\inf\{|1-p(g)|:g\in G\stm U\}=a >0.
$$
Let $f_U$ be a representation such that $p(g)=(f_U(g)v, v)$ for every $g\in G$
and some unit vector $v\in H$. Then the neighborhood 
$$
V=\{A\in U_s(H): |1-(Av, v)|<a\}
$$
of $1$ in $U_s(H)$ is as required.

It follows that the Hilbert direct sum of the representations $f_U$, 
taken over all neighborhoods $U$ of $e$, is topologically faithful. 
\end{proof}

Let $G$ be a locally compact Abelian group, $\hat G$ its Pontryagin dual.
Recall that continuous positive-definite functions on $\hat G$ are precisely
the Fourier transforms of positive measures on $G$ (Bochner's theorem). 
In particular, every positive-definite function $p:\R^k\to \C$ has the form
$$
p(x)=\int_{\R^k}\exp{(-2\pi i(x,y))}\,d\mu(y)
$$
for some positive measure $\mu$ on $\R^k$. 
Here $(x,y)=\sum_{i=1}^k x_iy_i$ for $x=(x_1,\dots,x_k)$ and  
$y=(y_1,\dots,y_k)$. Let $p(x)=\exp{(-\pi(x,x))}$.
Then $p$ is the Fourier transform of the Gaussian measure $p(y)dy$ and hence
is positive-definite. Thus for any $x_1,\dots,x_n\in \R^k$ the matrix
$(\exp{(-\norm{x_i-x_j}^2)})$ is positive. Here $\norm{\cdot}$ is the Hilbert
space norm on $\R^k$ defined by 
$\norm{y}=\sqrt{(y,y)}$.
If $H$ is an infinite-dimensional Hilbert space and $x_1,\dots,x_n\in H$,
then the finite-dimensional linear subspace of $H$ spanned by 
$x_1,\dots,x_n$ is isometric to $\R^k$ for some $k$. Thus we arrive at the 
following:

\begin{lemma}
\label{pos}
Let $H$ be a Hilbert space, and let $x_1,\dots,x_n$ be points in $H$.
Then the symmetric $n\ti n$-matrix $(\exp{(-\norm{x_i-x_j}^2)})$ is positive.
\end{lemma}

\section{Proof of the main theorem}
\label{3}

We prove a theorem which is more general than Theorem~\ref{main}.

\begin{thm}
\label{gen}
Let $(M,d)$ be a metric space, and let $G=\Is(M)$ be its group of isometries.
Suppose that there exists a real-valued positive-definite function $p:\R\to\R$
such that:
\begin{enumerate}
\item $p(0)=1$, and for every $\e>0$ we have $\sup\{p(x): |x|\ge \e\}<1$;
\item for every points $a_1,\dots,a_n\in M$ the symmetric real $n\ti n$-matrix
$(p(d(a_i,a_j)))$ is positive.
\end{enumerate}
Then the topological group $G$ is isomorphic to a subgroup of a unitary group.
\end{thm}

In virtue of Lemma~\ref{pos}, this theorem is indeed stronger than 
Theorem~\ref{main}: if $M=H$, 
we can take for $p$ the function defined by $p(x)=\exp{(-x^2)}$.

\begin{proof}
According to Proposition~\ref{crit}, we must show that there are sufficiently
many positive-definite functions on $G$.
For every $a\in M$ consider the function $p_a:G\to \R$ defined by
$p_a(g)=p(d(ga, a))$. This function is positive-definite. Indeed, 
let $g_1,\dots, g_n\in G$. Since each $g_i$ is an isometry of $M$, 
we have $d(g_i\obr g_ja,a)=d(g_ia, g_ja)$, and hence the matrix
$(p_a(g_i\obr g_j))=(p(d(g_i\obr g_ja,a)))=(p(d(g_ia,g_ja)))$ is positive
by our assumption. 

Let $U$ be a neighborhood of $1$ (= the identity map of $M$) in $G$. 
Without loss of generality we may assume that 
$$
U=\{g\in G: d(ga_i, a_i)<\e,\ i=1,\dots,n\}
$$
for some $a_1,\dots,a_n\in M$ and some $\e>0$. Let $t$ be the average of 
the functions $p_{a_1}, \dots, p_{a_n}$. Then $t$ is a positive-definite
function on $G$ such that 
$$
\sup\{t(g):g\in G\stm U\}<1.
$$
Indeed, let $\d>0$ be such that $p(x)\le 1-\d$ for every $x$ such that 
$|x|\ge\e$. If $g\in G\stm U$, then $d(ga_i,a_i)\ge\e$ for some $i$,
hence $p_{a_i}(g)=p(d(ga_i,a_i))\le 1-\d$ and 
$$
t(g)=\frac{1}{n}\sum_{k=1}^n p_{a_k}(g)\le 1-\frac{\d}{n}.
$$
We invoke Proposition~\ref{crit} to conclude that $G$ has a topologically
faithful representation.
\end{proof}

\end{document}